\documentclass[11pt]{article}

\usepackage{bbm}
\usepackage{amsmath}
\usepackage{amssymb}

\def\be{\begin{equation}}
\def\ee{\end{equation}}
\def\ba{\begin{eqnarray}}
\def\ea{\end{eqnarray}}
\def\lb{\label}
\def\nin{\noindent}

\def\dw{\!\!\downarrow}
\def\da#1{\!\!\downarrow\!\!^{\mbox{\rm\tiny ($#1$)}}}
\def\ua#1{\!\!\uparrow_{\!\mbox{\rm\tiny ($#1$)}}}
\def\lpl#1#2#3{{{#1}^+}^{#2}_{{#3}}}
\def\I{\ifmmode \mbox{\tiny I} \else \mbox{\footnotesize I}\fi}
\def\II{\ifmmode \mbox{\tiny II} \else \mbox{\footnotesize II}\fi}
\def\pha{\phantom{a}}

\newtheorem{proposition}{Proposition}
\newtheorem{corollary}[proposition]{Corollary}
\newtheorem{definition}[proposition]{Definition}

\def\teknum{\addtocounter{proposition}{1}\arabic{proposition}}

\input{amssym.def}
\input{amssym}

\makeatletter
\@addtoreset{equation}{section}
\makeatother

\begin{document}

\title{Bilinear identities on Schur symmetric functions}

\author{
\rule{0pt}{7mm} Dimitri Gurevich\thanks{gurevich@univ-valenciennes.fr}\\
{\small\it ISTV, Universit\'e de Valenciennes,
59313 Valenciennes, France}\\
\rule{0pt}{7mm} Pavel Pyatov\thanks{pyatov@theor.jinr.ru}\\
{\small\it Bogoliubov Laboratory of Theoretical Physics,
JINR, 141980 Dubna, Moscow region, Russia}\\
{\small \it and}\\
{\small\it Faculty of Mathematics, State University Higher School of Economics, 1001000 Moscow, Russia}\\
\rule{0pt}{7mm} Pavel Saponov\thanks{Pavel.Saponov@ihep.ru}\\
{\small\it Division of Theoretical Physics, IHEP, 142281
Protvino, Russia}
}
\date{}

\maketitle

\begin{abstract}
A series of bilinear identities on the Schur symmetric functions is
obtained with the use of Pl\"ucker relations.
\end{abstract}
\medskip

\leftline{{\bf Keywords:} Schur symmetric function, partition, Young diagram}
\smallskip

\leftline{2000 Mathematics Subject Classification: 05E05, 20B30} 

\section{Introduction}
\label{sec:1} 

In this note we present a wide class of bilinear
identities the Schur symmetric functions satisfy. The bilinear
identities are homogeneous second order polynomial relations with
integer coefficients, connecting different Schur functions. For the
detailed treatment of the Schur function theory, the corresponding
terminology, examples etc., see the monograph \cite{Mac}. Here we
give only a short list of definitions and key examples for
convenience of the reader.

A sequence of non-increasing non-negative integers
\[
\lambda = (\lambda_1,\lambda_2,\dots,\lambda_i,\dots),\quad
\lambda_1\ge\lambda_2\ge\dots\ge \lambda_i\ge \dots
\]
containing only finitely many non-zero terms is called {\it a
partition}. The total number of non-zero components,
$\ell(\lambda)$, is called {\it the height} of a given partition
$\lambda$
\[
\ell(\lambda) = n \;\Longleftrightarrow \;\lambda_n>0,\quad
\lambda_{n+1} = 0.
\]

Given a partition $\lambda$ with $\ell(\lambda) = n$, the Schur
symmetric function (actually, it is a polynomial)
$s_\lambda(t_1,\dots,t_m)$, where $m\ge \ell(\lambda)$, is an
element of the ring ${\Bbb Z}[t_1,\dots,t_m]$ defined as the ratio
of two determinants \cite{Mac}
\[
s_\lambda(t_1,\dots,t_m) = \left.\frac{\det\|t_i^{\lambda_j+m-j}\|}
{\det\|t_i^{m-j}\|}\right|_{1\le i,j\le m}.
\]
The set of Schur symmetric functions $s_\lambda(t_1,\dots,t_m)$
labeled by all partitions $\lambda$ with $\ell(\lambda)\le m$ forms
a ${\Bbb Z}$-basis of the subring of symmetric polynomials
\[
\Lambda_m = {\Bbb Z}[t_1,\dots,t_m]^{S_m}
\]
where the symmetric group $S_m$ acts on the polynomials from ${\Bbb
Z}[t_1,\dots,t_m]$ by the permutating the indeterminates.

The ring $\Lambda_m$ is graded
\[
\Lambda_m = \bigoplus_{k\ge 0}\Lambda_m^k,
\]
where $\Lambda_m^k$ consists of the homogeneous symmetric
polynomials of degree $k$. Then by a specific inverse limit (for
details, see \cite{Mac}) as $m\rightarrow \infty$ we pass from
$\Lambda_m^k$ to a graded ring $\Lambda$ called the \emph{ring of
symmetric functions in countably many indeterminates} $\{t_i\}_{i\in
\Bbb N}$. For each partition $\lambda$, the polynomials
$s_\lambda\in \Lambda_m$, define a unique element $s_\lambda\in
\Lambda$ called the \emph{Schur symmetric function in countably many
indeterminates}. Note that $s_\lambda\in \Lambda$ is no longer
polynomial (as well as other elements of the ring $\Lambda$). It is
a formal infinite sum of monomials, each of them being homogeneous
of degree $|\lambda| = \lambda_1+\dots+\lambda_n$. The Schur
symmetric functions form a ${\Bbb Z}$-basis of the ring $\Lambda$
and satisfy the Littlewood-Richardson multiplication rule \be
s_\lambda s_\mu = \sum_{\nu} C_{\lambda\mu}^{\nu}s_\nu, \label{LR}
\ee where the non-negative integers $C_{\lambda\mu}^{\nu}$ (the
Littlewood-Richardson coefficients) are calculated by some
combinatorial rule from partitions $\lambda$, $\mu$ and $\nu$.
Actually, the multiplication rule (\ref{LR}) can be taken 
for the formal definition of the ring $\Lambda$ in the ${\Bbb
Z}$-basis of Schur symmetric functions.

The bilinear identities we would like to discuss is another type of
relations among the Schur functions. As was mentioned at the
beginning of the section, they are of the form $p(\{s_{\lambda_i}\})
= 0$, where $p(\{x_i\})$ is a homogeneous second order polynomial (a
bilinear form) in its indeterminates with integer coefficients.
These identities follow, of course, from the multiplication rule
(\ref{LR}) but we use another technique to prove them.

As the first example of such identities we mention the bilinear
relations obtained in \cite{Kir}: \be s_{[m|n]}s_{[m|n]} =
s_{[m|n-1]}s_{[m|n+1]} + s_{[m-1|n]}s_{[m+1|n]}, \label{br-kir} \ee
where $[m|n]$ stands for the partition $(m^n)$ with $n$ components
equal to $m$. This identity connects the characters of the
irreducible representations of $SU(p+1)$, where
$s_{[m|n]}$ is a character of the $m$-th symmetric power of the
fundamental $SU(p+1)$ representation $\pi_n$ corresponding to
the signature $(1,1,\dots,1,0,\dots,0)$ ($n$ units, $1\le n\le
p$). The identity (\ref{br-kir}) played the key role in proving the
completeness of the Bethe vector set for the generalized
Heisenberg model. In the paper \cite{KR}, analogous bilinear
identities were obtained for the characters of symmetric powers of
fundamental representations of other classical Lie
groups (of B, C and D series).

In the work \cite{GPS} on quantum supermatrix algebras of $GL(m|n)$
type, we generalized the above identities to the products
$s_{[a|b]}\, s_{[m|n]}$ for arbitrary integers $1\le a\le m$ and
$1\le b\le n$: \be s_{[a|b]}\, s_{[m|n]}= \!\!\!\!\!\! \sum_{k=\atop
\max\{1,a+b-n\}}^a \!\!\!\!(-1)^{a-k} s_{[m|n]_{a+b-k}}\,
s_{[a-1|b-1]^{k-1}} + \!\!\!\!\!\! \sum_{k=\atop \max\{1, a+b-m\}}^b
\!\!\!\!(-1)^{b-k} s_{[m|n]^{a+b-k}}\, s_{[a-1|b-1]_{k-1}},
\label{gps-id} \ee where the symbols $[r|p]^k$ ($k\le r$) and
$[r|p]_k$ ($k\le p$) denote the partitions $((p+1)^k,p^{r-k})$
and $(p^r,k)$, respectively. These identities turned out to
be useful in studying the structure of the maximal commutative
subalgebras of the quantum supermatrix algebra.

In the work \cite{Kl}, identity (\ref{br-kir}) was generalized to
the product $s_\lambda s_\lambda$ for an arbitrary partition
$\lambda$. In the present paper, we give a different version of the
identity for the product $s_\lambda s_\lambda$. In contrast with the
result of \cite{Kl}, our formula admits the transposition of the
Young diagrams which parameterize the Schur functions. In other
words, given a bilinear identity for $s_\lambda s_\lambda$, we get a
true identity if we change all the partitions $\lambda$ by
their conjugates $\lambda'$ (see section \ref{sec:2} and
\cite{Mac}). In particular, if the Young diagram of the partition
$\lambda$ is symmetric under the transposition, the identity for
$s_\lambda s_\lambda$ is also symmetric.

M. Fulmek and M. Kleber have found the identities for the product of
two different Schur functions. Namely, in \cite{FK}, they proved
that \be
s_{(\lambda_1,\dots,\lambda_n)}s_{(\lambda_2,\dots,\lambda_{n+1})} =
s_{(\lambda_2,\dots,\lambda_n)}s_{(\lambda_1,\dots,\lambda_{n+1})}
+s_{(\lambda_2-1,\dots,\lambda_{n+1}-1)}s_{(\lambda_1+1,\dots,\lambda_n+1)},
\label{Dod} \ee where $(\lambda_1,\lambda_2,\dots,\lambda_{n+1})$ is
a partition, $n>0$ being an integer.

The series of the bilinear identities derived in this paper
considerably generalizes the identities (\ref{Dod}).

In the next section we introduce our notation and some key
operations with partitions. The third section is devoted to the
derivation of bilinear identities. The main results are formulated
in Proposition \ref{prop:1} and Corollary \ref{cor:f}.

\section{Definitions and notation}
\label{sec:2}

We use the terminology and definitions from the monograph
\cite{Mac}.

Let $\lambda = (\lambda_1,\dots,\lambda_n)$ be a partition of the
height $\ell(\lambda) = n$, that is $\lambda_n>0$. We omit the zero
components of $\lambda$. The Schur symmetric function corresponding
to the partition $\lambda$ can be expressed in terms of the complete
symmetric functions $h_k$ by means of the Jacobi-Trudi relations
\cite{Mac}: \be s_\lambda =\det \|h_{\lambda_i-i+j}\|_{1\le i,j\le
N}, \lb{JT} \ee where the index $i$ enumerates rows, the index $j$
enumerates columns, and $N\ge \ell(\lambda)=n$ is an arbitrary
positive integer. In the above formula it is assumed 
that $h_0\equiv 1$ and $h_k\equiv 0$ if $k<0$.
\medskip

\noindent{\bf Vectors $\mu$.} As is clear from the Jacobi-Trudi
determinant (\ref{JT}), any its row is completely defined by the
index of the first element of the row. Therefore, the Jacobi-Trudi
determinants and the corresponding Schur functions can be
unambiguously parameterized by the 
vectors $\mu \in {\Bbb Z}^N$ of the form \be
\mu=[\mu_1,\dots,\mu_N],\qquad \mu_i:= \lambda_i-i+1 \lb{def:mu} \ee
that is, $\mu =\lambda - \delta^{(N)}$, $ \delta^{(N)} =
[0,1,\dots,N-1]$. Unlike the partition $\lambda$, some of the
components of $\mu$ can be negative. Besides, the components of
$\mu$ form a {\it strictly} descending sequence
\[
\mu_1>\mu_2>\dots>\mu_N.
\]

To each partition $\lambda$ we assign its graphical image --- the
Young diagram (see \cite{Mac}). Below we denote the Young diagram of
the partition $\lambda$ by the same letter (when it does not lead to
a misunderstanding). Now we describe subsets of the Young diagram
$\lambda$ and define some operations with them; this will be used in
what follows.
\medskip

\nin{\bf The complete border strip.} Consider the Young diagram
corresponding to a partition $\lambda = (\lambda_1, \dots,
\lambda_n)$. Let us remove $\lambda_2-1$ boxes from the first row of
the diagram, starting from the first (the {\it left}-most) one. Then
we extend this procedure to the other rows removing
$\lambda_{k+1}-1$ boxes from the $k$-th row, $1\le k\le n-1$. We
leave the last $n$-th row unchanged.

This procedure results in a skew-diagram which will be referred to
as {\it the complete border strip}. Any nonempty proper subset of
the complete border strip will be called a {\it border strip}
provided this subset can be represented as the set-theoretical
difference $\lambda\setminus \nu$, where $\nu\subset \lambda$ is a
Young diagram completely contained in $\lambda$.

As an example, we consider the partition $(8,7,4^3,2^2)$. Its
Young diagram with the complete border strip marked by star signs
is depicted below:
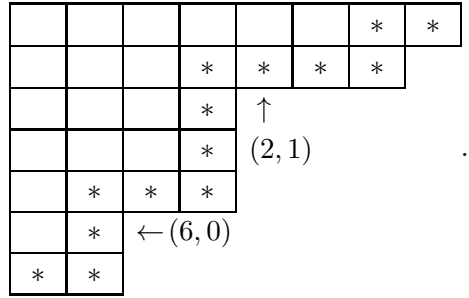
\[
\begin{array}{|c|c|c|c|c|c|c|c|}\hline
\hspace*{4mm}& \hspace*{4mm}&\hspace*{4mm} &\hspace*{4mm}
&\hspace*{4mm} &\hspace*{4mm} &\hspace*{1mm}*\hspace*{1mm}
&\hspace*{1mm}*\hspace*{1mm}\\ \hline
 & & &* &* &* &* & \multicolumn{1}{c}{} \\ \cline{1-7}
 & & &* &\multicolumn{1}{c}{ \uparrow}& \multicolumn{3}{c}{}\\
 \cline{1-4}
 & & &* &\multicolumn{3}{l}{(2,1)}\\ \cline{1-4}
 &* &* &* &\multicolumn{1}{c}{}\\ \cline{1-4}
 &* &\multicolumn{3}{l}{\leftarrow \!(6,0)}
 &\multicolumn{3}{c}{}\\ \cline{1-2}
 *&*\\ \cline{1-2}
\end{array}.
\]

We accept the following indexation of the boxes in the complete
border strip. As follows from the definition, in the $r$-th row of
the Young diagram $\lambda$, the boxes of the complete border strip
occupy positions from the $\lambda_{r+1}$-th column till the
$\lambda_r$-th one (counting from left to right). So, these boxes in
the $r$-th row can be enumerated by the number $s$ such that $0\le
s\le \lambda_r-\lambda_{r+1}$. A box of the complete border strip
situated in the $r$-th row and in the $(\lambda_{r+1}+s)$-th column
will be represented by an ordered pair of nonnegative integers
$(r,s)$. In the above example of the Young diagram, we show the
coordinate pairs of two boxes in the complete border strip.

\medskip

\nin{\bf The peeling.} Let us remove the complete border strip from
the Young diagram $\lambda$. The new diagram thus obtained will be
denoted by the symbol $\lambda\dw$. We say that $\lambda\dw$ is
obtained from $\lambda$ by peeling the complete border strip  off.
Note that the diagram $\lambda\dw$ can be the empty set if $\lambda$
is a simple hook diagram:
\[
(k,1^m)\dw = \emptyset\quad \text{for all~~}k,m\ge 0.
\]
It is not difficult to see that the diagram $\lambda\dw$ can be
obtained by removing the first row
and the first column from $\lambda$. As a consequence, the
height of $\lambda\dw$ is always
less than that of $\lambda$:
\[
\ell(\lambda\dw)\le \ell(\lambda)-1.
\]

Turning to the components of the partition $\lambda$, we get the
following structure of the partition $\lambda\dw$ \be \lambda =
(\lambda_1,\lambda_2,\dots,\lambda_n)\quad \rightarrow\quad
\lambda\dw = (\lambda_2-1,\lambda_3-1,\dots, \lambda_n-1,0).
\lb{srez} \ee The corresponding $\mu$-vectors (\ref{def:mu}) are
connected with each other by a simple transformation \be
\mu=[\mu_1,\mu_2,\dots,\mu_N]\quad\rightarrow\quad \mu\dw =
[\mu_2,\mu_3,\dots,\mu_N,-N+1]. \lb{ps} \ee In other words, the
components of $\mu$ are just shifted one position to the left, the
component $\mu_1$ disappears, and on the last place we get the
number $1- N$.

Consider now the peeling a border strip off, or a {\it partial}
peeling. In this case, we have to indicate the direction of the
peeling, that is we consider a \emph{partial up-peeling} and a
\emph{partial down-peeling}.

Let us fix a box $(r,s)$ in the complete border strip of a Young
diagram $\lambda$. Starting from the box $(r,s)$, we remove all the
boxes of the complete border strip lying to the left and down of the
chosen box. That is, we remove all the boxes $(r,t)$ with $0\le t\le
s$ and $(p,t)$ with $p>r$. This procedure will be called the {\it
partial down-peeling} from the starting box $(r,s)$. We will only be
interested in down-peelings that transform a Young diagram to a
Young diagram. For this to be true, the starting box $(r,s)$ of the
partial down-peeling must be the {\it right-most} box in the $r$-th
row. In other words, the number $s$ must take the maximal possible
value $s=\lambda_{r}-\lambda_{r+1}$. To simplify the expressions, we
omit this $s$ in notation and denote the diagram (and the partition)
obtained from the diagram $\lambda$ by the partial down-peeling from
the box $(r,\lambda_{r}-\lambda_{r+1})$ by the symbol
$\lambda\da{r}$. The components of the partition $\lambda\da{r}$
read \be \lambda\da{r} =
(\lambda_1,\dots,\lambda_{r-1},\lambda_{r+1}-1,
\dots,\lambda_{n}-1,0), \label{l-p-sr} \ee while for the components
of the corresponding $\mu$-vector $\mu\da{r}$ we obtain \be
\mu\da{r} = [\mu_1,\dots,\mu_{r-1},\mu_{r+1},\dots,\mu_N,-N+1].
\lb{ch-sr-vn} \ee Same as the peeling the complete border strip off,
the partial down-peeling decreases the height of the diagram at
least by one: $\ell(\lambda\da{r})\le \ell(\lambda) -1$.

The partial up-peeling is defined in an analogous way. We fix a
starting box $(r,s)$ in the complete border strip of a diagram
$\lambda$ and remove all the boxes $(r,t)$ with $t\ge s$ and $(p,t)$
with $p<r$. That is we remove all the boxes of the complete border
strip, lying to the right and up of the chosen starting box. This
procedure will be called the {\it partial up-peeling} from the
starting box $(r,s)$. In what follows we will be interested only in
partial up-peelings that do not destroy the structure of Young
diagrams. Therefore, the starting box $(r,s)$ of the up-peeling must
be chosen in such a way that there are no box of the diagram
directly {\it under} it. This is only possible if $\lambda_{r}>
\lambda_{r+1}$ and, besides, $s\ge 1$. The Young diagram (and the
partition) obtained from the diagram $\lambda$ by the partial
up-peeling from the starting box $(r,s)$ will be denoted by the
symbol $\lambda\ua{r,s}$.

The component structure of the partition $\lambda\ua{r,s}$ is as
follows
\be
\lambda\ua{r,s} =
\begin{array}[t]{ccccccc}
(\lambda_2-1, & \dots, & \lambda_{r}-1, &
\lambda_{r+1}+s-1, & \lambda_{r+1}, & \dots , & \lambda_n)\\
\mbox{\footnotesize 1} & & \mbox{\footnotesize $r-1$}
 & \mbox{\footnotesize $r$} & \mbox{\footnotesize $r+1$} &
 & \mbox{\footnotesize $n$}
\end{array},\quad
1\le s \le \lambda_{r-1}-\lambda_r, \label{l-up-sr} \ee where in the
second line we have written the ordinal numbers of the corresponding
components to clarify the structure. For the corresponding vector
$\mu$, we get the following expression \be \mu\ua{r,s} =
\begin{array}[t]{ccccccc}
[\mu_2, & \dots, & \mu_{r}, &
\mu_{r+1}+s, & \mu_{r+1}, & \dots , & \mu_N]\\
\mbox{\footnotesize 1} & & \mbox{\footnotesize $r-1$}
 & \mbox{\footnotesize $r$} & \mbox{\footnotesize $r+1$} &
 & \mbox{\footnotesize $N$}
 \end{array},\quad 1\le s \le \mu_{r-1}-\mu_r-1,
 \label{ch-sr-up}
\ee

\medskip

\nin{\bf Adding a border strip to diagram.} Consider the Young
diagram, corresponding to a partition
$\lambda=(\lambda_1,\dots,\lambda_n)$. Choose $m\le n-1$ consecutive
rows with numbers $r,r+1,\dots , r+m-1$, where $2\le r\le n-m +1$.
We are going to add boxes in the chosen rows in such a way that the
result would be a Young diagram, and, besides, the added boxes would
form a connected border strip in the new diagram. The restriction on
the number of rows means that we do not add boxes into the first
line of $\lambda$ ($r\ge 2)$ and that we do not increase the height
of the diagram ($r\le n-m+1$). Below we use the shorthand notation
$r_m:=r+m-1$.

It turns out to be convenient to treat the first (the  left-most)
box added into the $r_m$-th row as {\it the beginning} (or the
first) box of the strip.

The last (the  right-most) box added into the $r$-th row will be
treated as {\it the end} (or the last) box of the strip. The
beginning of the added strip can be placed in any row of $\lambda$
(except for the above restriction on number) with the only
requirement that the first added box must appear in the
$(\lambda_{r_m}+1)$-th column (to preserve the correct structure of
the Young diagram). As for the end of the strip, it can be situated
only in the row which is shorter than its preceding row:
$\lambda_r<\lambda_{r-1}$.

\label{str-pol} The number of boxes added into the $(r+i)$-th row
reads as follows \be p_i = \lambda_{r+i-1} - \lambda_{r+i}+1,\quad
1\le i\le m-1. \lb{qi} \ee Into the last, $r$-th, row we add $p_0 =
t$ boxes, where $1\le t\le \lambda_{r-1}-\lambda_r$. Therefore, the
total amount of boxes added is equal to
\[
p =\sum_{i=0}^{m-1}p_i = \lambda_r - \lambda_{r_m} +t+m-1
=\mu_r-\mu_{r_m}+t.
\]
Here is an example of adding a border strip for the case
$\lambda=(8,7,4^3,2^3)$, $r=3$, $m=5$ and $t=2$:
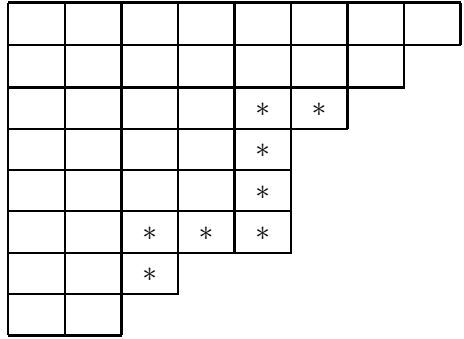
\[
\begin{array}{|c|c|c|c|c|c|c|c|}\hline
\hspace*{4mm} & \hspace*{4mm}&\hspace*{4mm} &\hspace*{4mm}
&\hspace*{4mm} &\hspace*{4mm} &\hspace*{4mm} &\hspace*{4mm}\\ \hline
 & & & & & & & \multicolumn{1}{c}{} \\ \cline{1-7}
 & & & &\multicolumn{1}{c|}{*}&\multicolumn{1}{c|}{*}\\
 \cline{1-6}
 & & & &\multicolumn{1}{c|}{*}\\ \cline{1-5}
 & & & &\multicolumn{1}{c|}{*}\\ \cline{1-5}
 & &*&* &*\\ \cline{1-5}
 & &\multicolumn{1}{c|}{*} \\ \cline{1-3}
 &\\ \cline{1-2}
\end{array}
\]
Here stars denote the added boxes.

The symbol $\lpl{\lambda}{t}{(r, m)}$ will stand for the diagram
(and the partition) obtained from the diagram $\lambda$ by adding a
border strip of $m$ rows from $r$ to $r_m=r+m-1$ with $t$ boxes in
the end row $r$. If we add several (say $k$) disconnected border
strips, the notation is obviously generalized to
$\lpl{\lambda}{t_1\quad\dots\quad t_k}{(r_1,m_1)\dots (r_k,m_k)}$.

The components of the partition $\lpl{\lambda}{t}{(r,m)}$ read
(recall that $1\le t\le \lambda_{r-1}-\lambda_r$)
\be
\begin{array}[t]{cccccccccccc}
\lpl{\lambda}{t}{(r,m)}& = &
(\lambda_1,\dots,\lambda_{r-1},&\lambda_r+t,&
\lambda_r+1,&\lambda_{r+1}+1,& \dots, &
\lambda_{r_m-1}+1,&\lambda_{r_m+1},&
\dots,&\lambda_n)\\
& & &\mbox{\footnotesize $r$} &
  \mbox{\footnotesize $r+1$} &
  \mbox{\footnotesize $r+2$} &\dots &
  \mbox{\footnotesize $r_m$}&
\mbox{\footnotesize $r_m+1$} &\dots
 &\mbox{\footnotesize $n$}.
\end{array}
\lb{la-pol} \ee Here in the second line we have written the ordinal
numbers of the corresponding  components.

The component structure of the corresponding vector
$\lpl{\mu}{t}{(r,m)}$ is more transparent \be
\begin{array}[t]{cccccccccccc}
\lpl{\mu}{t}{(r,m)} & = &
[\mu_1,\dots,\mu_{r-1},& \mu_r+t,&
\mu_r,& \mu_{r+1},& \dots, &
\mu_{r_m-1},& \mu_{r_m+1},&
\dots,& \mu_N]\\
 & & &\mbox{\footnotesize $r$} &
  \mbox{\footnotesize $r+1$} &
  \mbox{\footnotesize $r+2$} & \dots &
   \mbox{\footnotesize $r_m$}&
 \mbox{\footnotesize $r_m+1$} & \dots
 & \mbox{\footnotesize $N$}
\end{array}.
\lb{mu-pol} \ee As we see, the changes take place only for the
components from $\mu_r$ to $\mu_{r_m}$. Namely, the string of
components $\mu_r,\dots,\mu_{r_m-1}$ shifts one position to the
right, in the $r$-th place (the end row of the added strip) we get
the new component $\mu_{r}+t$ and the component $\mu_{r_m}$ (the
beginning row of the strip) disappears.

\section{Bilinear identities}
\label{sec:3}

The bilinear identities on the Schur symmetric functions follow from
the Jacobi-Trudi determinant formula (\ref{JT}) and the Pl\"ucker
relation on the product of two determinants (for details, see
\cite{Sturm}). Let us formulate the corresponding statement for the
reader's convenience.

Consider a pair of $\, p\times p\, $ matrices $A =
\|a_{ij}\|_{i,j=1}^p$ and $B=\|b_{ij}\|_{i,j=1}^p$. Let $a_{i*}$
denote the $i$-th row of the matrix $A$. Introduce the following
notation: \be \lb{notat} \det A\, :=\, |A|\, , \qquad A\, :=\,
\left(
\begin{array}{ccccc}
a_{1*} &\dots &a_{i*} & \dots & a_{p*}\\
1 &\dots &i & \dots & p
\end{array}
\right), \ee where the last symbol contains a detailed information
on the row content of $A$. Namely, it says that the row $a_{i*}$ is
located in the $i$-th place in the matrix $A$ (when counting from
the top down).

Let us fix a set of integer data $\{k\,|\,r_1,r_2,\dots ,r_k\}$,
where $1\leq k\leq p$ and $1\leq r_1<\dots < r_k\leq p$. Given these
data, the Pl{\"u}cker relation reads
\begin{eqnarray}
|A| |B| = \sum_{1\le s_1<\dots <s_k\le p}&&
\left|
\begin{array}{ccccccccc}
a_{1*} & \dots & b_{s_1*} & \dots & b_{s_2*} &
\dots
& b_{s_k*} &\dots & a_{p*} \\
1 & \dots & r_1 & \dots & r_2 & \dots & r_k & \dots
& p
\end{array}
\right|\times \nonumber\\
&&\left|
\begin{array}{ccccccccc}
b_{1*} & \dots & a_{r_1*} & \dots & a_{r_2*} &
\dots
& a_{r_k*} &\dots & b_{p*} \\
1 & \dots & s_1 & \dots & s_2 & \dots & s_k & \dots
& p
\end{array}
\right| ,
\label{Pluk}
\end{eqnarray}
where the sum is taken over all possible sets $\{k\,|\,s_1,\dots ,s_k\}$.

Now we can obtain a bilinear identity, connecting the Schur
symmetric functions labeled by a partition $\lambda$ and the
partition $\lpl{\lambda}{t_1\quad\dots\quad t_k}{(r_1,m_1)\dots
(r_k,m_k)}$. Here we assume that the structure of the diagram
$\lambda$ allows adding $k$ border strips of the indicated size and
location.

\begin{proposition}
\label{prop:1} In the Young diagram corresponding to a partition
$\lambda = (\lambda_1,\dots,\lambda_n)$, let there exist $k\ge 1$
rows with numbers $2\le r_1 < r_2 < \dots < r_k \le r_{k+1}:=n$
possessing the property
\[
\lambda_{r_i} < \lambda_{r_i-1},\quad 1\le i \le k.
\]
Let the integers $t_i,m_i$, where $1\le i\le k$, satisfy the
restrictions
\[
1\le t_i\le \lambda_{r_i-1}-\lambda_{r_i},\quad 1\le m_i\le r_{i+1} -r_i,\quad 1\le i \le k.
\]
Then the Young diagram $\lpl{\lambda}{t_1\quad\dots\quad
t_k}{(r_1,m_1)\dots (r_k,m_k)}$ can be defined and the following
bilinear identity on the Schur symmetric functions holds {\rm \be
s_{\raisebox{-2pt}{$_\lambda$}}s_{\lpl{\lambda}{t_1\quad\dots\quad
t_k}{(r_1,m_1) \dots (r_k,m_k)}\,\,\dw} =
s_{\lpl{\lambda}{t_1\quad\dots\quad t_k}{(r_1,m_1) \dots (r_k,m_k)}}
s_{\raisebox{-2pt}{$_{\lambda\,\,\dw}$}} +\sum_{p=1}^k
s_{\lpl{\lambda}{t_1\quad\dots\quad t_k}{(r_1,m_1) \dots
(r_k,m_k)}\,\,\downarrow^{\mbox{\tiny$(r_p)$}}}s_{\raisebox{-2pt}
{$_{\lambda\,\,\ua{r_p-1,t_p}}$}}. \label{main-id} \ee }
\end{proposition}

\nin{\bf Proof.} To prove the proposition we use the Jacobi-Trudi
formulae for the Schur functions and the Pl\"ucker relation for the
product of two determinants. In so doing, we shall parameterize the
rows of the Jacobi-Trudi determinants in (\ref{Pluk}) by components
of vectors $\mu$ defined in (\ref{def:mu}).

First of all, we inspect the structure of the Jacobi-Trudi
determinants in the left-hand side of (\ref{main-id}) in order to
find the set of rows to be exchanged in accordance with the
Pl\"ucker relation. Taking into account expression (\ref{mu-pol})
for the $\mu$-vector of the diagram with added border strip and
expression (\ref{ps}) for the peeling the complete border strip off,
we have
\begin{eqnarray*}
&s_{\raisebox{-2pt}{$_\lambda$}}s_{\lpl{\lambda}{t_1\quad\dots\quad t_k}{(r_1,m_1)
\dots (r_k,m_k)}\,\,\dw} =\left|\begin{array}{cccccccc}
\mu_1 &\dots & \mu_{r_i}&\mu_{r_i+1} &\dots&\mu_{r_{m_i}}&
\dots&\mu_N\\
1& \dots & r_i &r_i+1& \dots & r_{m_i} & \dots & N
\end{array}\right|_{\,\mbox{\sf \footnotesize I}}\times&
\\
&
\left|\begin{array}{cccccccccccc}
\mu_2 & \dots & \mu_{r_i-1}& \mu_{r_i}+t_i&
\mu_{r_i}& \mu_{r_i+1}&\dots & \mu_{r_{m_i}-1}&\mu_{r_{m_i}+1}&
\dots&\mu_N& -N+1\\
1 & \dots & r_i-2 & r_i-1 & r_i &r_i+1 & \dots & r_{m_i}-1 & r_{m_i}&
\dots & N-1& N
\end{array}\right|_{\,\mbox{\sf \footnotesize II}},&
\end{eqnarray*}
where we explicitly indicated the components containing the $i$-th
part of the added border strip. Recall that it is located in rows
between $r_i$ and $r_{m_i} = r_i+m_i-1$. The indices {\sf I} and
{\sf II} were introduced for convenience of references.

Let us take the data $\{k\,|\,r_{m_1},r_{m_2},\dots r_{m_k}\}$ to
indicate the $k$ rows of the determinant {\sf I} to be exchanged
with all possible sets of $k$ rows of the determinant {\sf II} in
accordance with the Pl\"ucker relation (\ref{Pluk}). It is not
difficult to see that in the right hand side of the Pl\"ucker
relation applied to the above product of the determinants {\sf I}
and {\sf II} there are only $(k+1)$ nonzero terms. They correspond
to the exchange of the rows $r_{m_1},r_{m_2},\dots, r_{m_k}$ of the
determinant {\sf I} with rows $r_1-1,r_2-1,\dots, r_k-1$ and $N$ of
the determinant {\sf II}. The other terms vanish since the
determinants obtained in exchanging procedure possess at least two
identical rows.

The nonzero terms correspond to the following ways of row exchange
\[
\begin{array}{l@{\hspace*{5mm}}l}
\left[ \begin{array}{c}
\mu_{r_{m_{i}}}\\
r_{m_i}
\end{array} \right]_{\mbox{\sf \footnotesize I}}
\;\longleftrightarrow\; \left[ \begin{array}{c}
\mu_{r_i}+t_i\\
r_i-1
\end{array} \right]_{\mbox{\sf \footnotesize II}},
\quad 1\le i\le
k&\mbox{ placement A}\\
&\\
\hspace*{-5pt}\left.\begin{array}{lcc@{\quad}l}
\left[ \begin{array}{c}
\mu_{r_{m_{i}}}\\
r_{m_i}
\end{array} \right]_{\mbox{\sf \footnotesize I}}
&\longleftrightarrow& \left[ \begin{array}{c}
\mu_{r_i}+t_i\\
r_i-1
\end{array} \right]_{\mbox{\sf \footnotesize II}},
& 1\le i\le p-1 \\
\rule{0pt}{9mm}
\left[ \begin{array}{c}
\mu_{r_{m_{j}}}\\
r_{m_j}
\end{array} \right]_{\mbox{\sf \footnotesize I}}
&\longleftrightarrow& \left[ \begin{array}{c}
\mu_{r_{j+1}}+t_{j+1}\\
r_{j+1}-1
\end{array} \right]_{\mbox{\sf \footnotesize II}},
& p\le j\le k-1 \\
\rule{0pt}{9mm}
\left[ \begin{array}{c}
\mu_{r_{m_{k}}}\\
r_{m_k}
\end{array} \right]_{\mbox{\sf \footnotesize I}}
&\longleftrightarrow& \left[ \begin{array}{c}
-N+1\\
N
\end{array} \right]_{\mbox{\sf \footnotesize II}}
\end{array}\right\}
&\mbox{placements $B_p$, $1\le p\le k$}
\end{array}
\]
The row exchange in accordance with the placement A gives the first
term in the right hand side of (\ref{main-id}). Indeed, after such
an exchange the typical part of the determinant {\sf I} takes the
form
\[
\left|\begin{array}{cccccc}
\dots & \mu_{r_i}&\mu_{r_i+1} &\dots&\mu_{r_i}+t_i& \dots\\
\dots & r_i &r_i+1& \dots & r_{m_i} & \dots
\end{array}\right|_{\,\mbox{\sf \footnotesize I}}.
\]
Now we have to make the cyclic permutation of rows placing the
component $\mu_{r_i}+t_i$ to the $r_i$-th row. This gives the sign
factor $(-1)^{r_{m_i}-r_i} = (-1)^{m_i-1}$ and, according to
(\ref{mu-pol}), the structure of the determinant {\sf I} corresponds
to the Schur function $s_{\lpl{\lambda}{t_1\quad\dots\quad
t_k}{(r_1,m_1) \dots (r_k,m_k)}}$. As for the typical part of the
second determinant, we get after the row exchange
\[
\left|\begin{array}{ccccccccccc}
\mu_2 & \dots & \mu_{r_{m_i}}&
\mu_{r_i}& \dots & \mu_{r_{m_i}-1}&\mu_{r_{m_i}+1}& \dots&\mu_N& -N+1\\
1 & \dots & r_i-1 & r_i & \dots & r_{m_i}-1&r_{m_i} & \dots & N-1& N
\end{array}\right|_{\,\mbox{\sf \footnotesize II}}.
\]
Here we also have to make the cyclic permutation of rows from
$(r_i-1)$ to $(r_{m_i}-1)$ placing the component $\mu_{r_{m_i}}$ to
the $(r_{m_i}-1)$-th row. This generates the sign factor
$(-1)^{m_i-1}$ which compensates the same factor of the determinant
{\sf I}. As for the structure of the determinant {\sf II}, it
corresponds to $s_{\lambda\,\,\dw}$ as directly follows from
(\ref{ps}).

Turn now to a placement of $B_p$ type for some fixed integer $p$
such that $1\le p\le k$. We first consider the changes in the
determinant {\sf I}. The rows $r_{m_1}$ to $r_{m_{p-1}}$ are
exchanged in the same way as in the placement A giving rise to the
following typical parts corresponding to added border strips
\[
\left|\begin{array}{ccccccc}
\mu_1&\dots & \mu_{r_i}+t_i&\mu_{r_i} &\dots&\mu_{r_{m_i}-1}& \dots\\
1&\dots & r_i &r_i+1& \dots & r_{m_i} & \dots
\end{array}\right|_{\,\mbox{\sf \footnotesize I}},\quad 1\le i\le p-1,
\]
with the sign factor $(-1)^{m_i-1}$ for each strip lying in rows
$r_i$ to $r_{m_i}$. The remaining part of the determinant {\sf I}
can be transformed to
\[
\left|\begin{array}{ccccccccccccc} \dots&\mu_{r_{m_p}-1}
&\mu_{r_{m_p}+1}& \dots&\mu_{r_j}+t_j&\mu_{r_j}
&\dots&\mu_{r_{m_j}-1}& \dots
&\mu_{r_{m_k}-1}&\mu_{r_{m_k}+1}&\dots&-N+1\\
\dots&r_{m_p}-1 &r_{m_p} &\dots & r_j-1 &r_j& \dots & r_{m_j}-1 &
\dots &r_{m_k}-1&r_{m_k}&\dots&N
\end{array}\right|_{\,\mbox{\sf \footnotesize I}}
\]
with the sign factors $(-1)^{r_{j+1}-r_{m_j}-1}$ and $p\le j\le k-1$
which originate from the cyclic permutation of rows from $r_{m_j}$
till $(r_{j+1}-1)$. This permutation results in moving the component
$\mu_{r_{j+1}}+t_{j+1}$ from the $r_{m_j}$-th row to the
$(r_{j+1}-1)$-th one. We have also a sign factor $(-1)^{N-r_{m_k}}$
since the component $(-N+1)$ moved from the $r_{m_k}$-th row to the
last, $N$-th, row. Finally, taking into account the structure of the
partial down-peeling (\ref{ch-sr-vn}), we see that, up to the above
sign factors, the determinant ${\sf I}$ represents the following
Schur symmetric function
\[
\left|\begin{array}{cccccccc}
\mu_1 &\dots & \mu_{r_i}&\mu_{r_i+1} &\dots&\mu_{r_{m_i}}&
\dots&\mu_N\\
1& \dots & r_i &r_i+1& \dots & r_{m_i} & \dots & N
\end{array}\right|_{\,\mbox{\sf \footnotesize I}}
\stackrel{B_p}{\longrightarrow} s_{\lpl{\lambda}{t_1\quad\dots\quad
t_k}{(r_1,m_1) \dots
(r_k,m_k)}\,\,\downarrow^{\mbox{\tiny$(r_p)$}}}.
\]

Consider now the changes in the determinant {\sf II} under the row
exchange of the same $B_p$ type. The part of the determinant
containing the rows $r_i-1$ for $1\le i\le p-1$ can be expressed in
the following form
\[
\left|\begin{array}{cccccccc}
\mu_2 &\dots & \mu_{r_{m_i}}&\mu_{r_i} &\dots&\mu_{r_{m_i}-1}&\mu_{r_{m_i}+1}&
\dots\\
1& \dots & r_i-1 &r_i& \dots & r_{m_i}-1& r_{m_i}& \dots
\end{array}\right|_{\,\mbox{\sf \footnotesize II}}.
\]
Here we have to rearrange the rows from $(r_i-1)$ till $(r_{m_i}-1)$
by cyclic permutation in order to move the component $\mu_{r_{m_i}}$
to the $(r_{m_i}-1)$-th row. This gives rise to the sign factor
$(-1)^{r_i-r_{m_i}} = (-1)^{m_i-1}$ for each $1\le i\le p-1$. The
sign factors compensate the analogous sign factors originated from
the determinant {\sf I}.

The rest part of the determinant {\sf II} reads ($p\le j\le k-1$)
\[
\left|\begin{array}{ccccccccccc}
\dots& \mu_{r_p}+t_p&\mu_{r_p} &\dots & \mu_{r_{m_j}}&\mu_{r_{j+1}+1}&\dots &
\mu_{r_{m_{j+1}}-1}&\mu_{r_{m_{j+1}}+1}&\dots&\mu_{r_{m_k}}\\
\dots & r_p-1 &r_p & \dots & r_{j+1}-1& r_{j+1}& \dots &r_{m_{j+1}}-1&r_{m_{j+1}}&\dots
&N
\end{array}\right|_{\,\mbox{\sf \footnotesize II}}.
\]
On moving the component $\mu_{r_{m_j}}$ from the $(r_{j+1}-1)$-th
row to the $r_{m_j}$-th one we get the sign factor
$(-1)^{r_{j+1}-r_{m_{j}}-1}$ for each $p\le j\le k-1$. Also we get
the factor $(-1)^{N-r_{m_k}}$ since the component $\mu_{r_{m_k}}$
moved from the last, $N$-th, row to the row $r_{m_k}$. All these
sign factors exactly compensate the corresponding sign factors
appearing in the determinant {\sf I}. The final structure of the
determinant {\sf II} is as follows:
\[
\left|\begin{array}{ccccccccccc} \mu_2 &\dots & \mu_{r_i}&\dots
&\mu_{r_{p-1}}&\mu_{r_p}+t_p&\mu_{r_p}&\dots&\mu_{r_j}&\dots& \mu_N\\
1& \dots & r_i-1 &\dots & r_p-2&r_p-1& r_p& \dots& r_j&\dots& N
\end{array}\right|_{\,\mbox{\sf \footnotesize II}}.
\]
On comparing the above determinant with (\ref{ch-sr-up}), we
conclude that under the row exchange of the $B_p$ type the
determinant {\sf II} transforms to the Schur function
$s_{\raisebox{-2pt}{$_{\lambda\,\,\ua{r_p-1,t_p}}$}}$ (up to the
sign factors compensated by the corresponding factors of the
determinant {\sf I}).

At last, summing over all placements of the $B_p$ type and adding
the result of the placement $A$ we get the final formula
(\ref{main-id}). \hfill\rule{6.5pt}{6.5pt}
\medskip

Consider now some important corollaries of Proposition \ref{prop:1}.

\begin{corollary}{\rm
 \label{cor:1}
The identity (\ref{main-id}) is preserved under the simultaneous
transposition of all the Young diagrams parameterizing the Schur
functions in (\ref{main-id}). }
\end{corollary}

\nin{\bf Proof.} Recall (see \cite{Mac}) that the partition
$\lambda'$ is said to be the {\it conjugate} of a given partition
$\lambda$ if the Young diagram $\lambda'$ is obtained from the Young
diagram $\lambda$ by the transposition with respect to the main
diagonal. In other words,
\[
\lambda'_i = \#(j\mid \lambda_j\ge i).
\]
The key point in the proof of the Corollary \ref{cor:1} is the
following Jacobi-Trudi determinant formula for the Schur symmetric
function $s_\lambda$ \be s_\lambda =
\det\|e_{\lambda'_i-i+j}\|_{1\le i,j,\le M}, \label{JT-el} \ee where
$e_k$ is the $k$-th elementary symmetric function, and $M\ge
\ell(\lambda') = \lambda_1$ is an arbitrary positive integer. Here,
as well as in relation (\ref{JT}), we set: $e_k\equiv 0$ for $k<0$
and $e_0\equiv 1$.

The proof of Proposition \ref{prop:1} is based on formula
(\ref{JT}), which contains the complete symmetric functions $h_k$.
But we do not use any specific properties of these functions in
course of the proof. The functions $h_k$ are just the matrix
elements of determinants in the Pl\"ucker relation. If we change all
the complete symmetric functions $h_{\lambda_i-i+j}$ for
$e_{\lambda_i-i+j}$ the identity (\ref{main-id}) still remains true
{\it determinant} identity. The interpretation of the determinants
involved will, however, be different. As can be seen from
(\ref{JT-el}), the determinants will now parameterize the Schur
functions corresponding to the conjugate partitions $\lambda'$.
\hfill\rule{6.5pt}{6.5pt}
\medskip

Another useful consequence of the proof of Proposition \ref{prop:1}
is a possibility to remove the first line or the first column of
some partitions and get a new identity. Indeed, as can be easily
seen from the proof, the first row of the Jacobi-Trudi determinant
corresponding to the Schur function $s_\lambda$ (the component
$\mu_1$) does not play an active role in the calculations. In
principle, it can be changed for an arbitrary row and identity
(\ref{main-id}) will be still valid as the {\it determinant}
identity (though the interpretation of the corresponding
determinants as Schur functions will be lost in general). But if we
change the row $\mu_1$ by the $N$-dimensional row
$(1,0,\dots,0)$, the determinants $s_\lambda$, $s_{\lambda^+}$ and
$s_{\lambda^+\!\downarrow^{(r)}}$ can be interpreted as the Schur
functions corresponding to the partition with the first component
removed. Here is an example of the procedure:
\[
s_{(\lambda_1,\dots,\lambda_n)}=\det\|h_{\lambda_i-i+j}\|
\quad\longrightarrow\quad
\left|\begin{array}{llll}
1&0&\dots&0\\
h_{\lambda_2-1}&h_{\lambda_2}&\dots&h_{\lambda_2+n-2}\\
\vdots&\dots&\vdots&\vdots\\
h_{\lambda_n-n+1}&h_{\lambda_n-n+2}&\dots&h_{\lambda_n}\\
\end{array}\right| =s_{(\lambda_2,\dots,\lambda_n)}.
\]
Due to Corollary~\ref{cor:1} the same is true for removing the first column
in the diagram $\lambda$. Therefore, the following corollary holds true.

\begin{corollary}
\label{cor:2}{\rm Let $\lambda=(\lambda_1,\dots,\lambda_n)$ be a
partition satisfying the conditions of Proposition~\ref{prop:1}.
Denote by $\bar\lambda$ the partition obtained from $\lambda$ by
removing the first line or the first column from the Young diagram
$\lambda$, that is
\[
\bar\lambda = (\lambda_2,\dots,\lambda_n)\quad{\rm or}\quad
\bar\lambda = (\lambda_1-1,\lambda_2-1,\dots,\lambda_n-1).
\]
Then identity (\ref{main-id}) implies that
\def\tld{\bar\lambda}
\be
s_{\raisebox{-1pt}{$_{\bar\lambda}$}}s_{\lpl{\lambda}{t_1\quad\dots\quad
t_k}{(r_1,m_1) \dots (r_k,m_k)}\,\,\dw} =
s_{{\overline{\lambda^+}}^{\;t_1\quad\dots\quad
t_k}_{\;\;(r_1,m_1)\dots (r_k,m_k)}}
s_{\raisebox{-2pt}{$_{\lambda\,\,\dw}$}} +\sum_{p=1}^k
s_{{\overline{\lambda^+}}^{\;t_1\quad\dots\quad
t_k}_{\;\;(r_1,m_1)\dots (r_k,m_k)}
\downarrow^{\mbox{\tiny$(r_p)$}}}
\,s_{\raisebox{-2pt}{$_{\lambda\,\,\ua{r_p-1,t_p}}$}}.
\label{cor-id} \ee Here ${\overline{\lambda^+}}$ and
${\overline{\lambda^+}}\!\! \downarrow^{(r_p)}$ are the Young
diagrams obtained from the diagrams $\lambda^+$ and
$\lambda^+\!\!\!\downarrow^{(r_p)}$ by removing the first row
(column). }
\end{corollary}

We give two examples illustrating the above formulae.\par
\medskip

\nin{\bf Example \teknum} Let $\lambda=(2,1,1)$, $k=1$, $r_1=2$,
$m_1=1$ in accordance with the notation of Proposition~\ref{prop:1}. That
is we add a single box in the second row of the Young diagram
$\lambda$. Then the main identity (\ref{main-id}) reads:
\[
s_{(2,1,1)}s_{(1)} = s_{(2,2,1)} + s_{(2)}s_{(1,1,1)}
\]
or, loosely denoting the Schur functions $s_\lambda$ by the
corresponding Young diagrams $\lambda$ (for more visual clarity)
\[
\begin{tabular}{|c|c|}\hline
 & \\ \hline
&\multicolumn{1}{c}{}\\ \cline{1-1}
&\multicolumn{1}{c}{}\\ \cline{1-1}
\end{tabular}\;\times\;
\begin{tabular}{|c|}\hline
  \\ \hline
\end{tabular}\; =\;
 \begin{tabular}{|c|c|}\hline
  & \\ \hline
  & \\ \hline
&\multicolumn{1}{c}{}\\ \cline{1-1}
\end{tabular}\quad+\quad
\begin{tabular}{|c|c|}\hline
  & \\ \hline
\end{tabular}\;\times\;
\begin{tabular}{|c|}\hline
  \\ \hline
  \\ \hline
  \\ \hline
\end{tabular}\;.
\]
On removing the first row ($\lambda\rightarrow \bar\lambda=(1,1)$), we get
\[
s_{(1,1)}s_{(1)} = s_{(2,1)} +s_{(1,1,1)},
\]
or, in the graphic form
\[
\begin{tabular}{|c|}\hline
  \\ \hline
  \\ \hline
\end{tabular}\;\times \;
\begin{tabular}{|c|}\hline
  \\ \hline
\end{tabular}\;=\;
\begin{tabular}{|c|c|}\hline
  & \\ \hline
  &\multicolumn{1}{c}{}\\ \cline{1-1}
\end{tabular}\;+\;
\begin{tabular}{|c|}\hline
  \\ \hline
  \\ \hline
  \\ \hline
\end{tabular}\;.
\]
On removing the first column, we find
\[
s_{(1)}s_{(1)} = s_{(1,1)} + s_{(2)},
\]
or, in the graphic form
\[
\begin{tabular}{|c|}\hline
  \\ \hline
\end{tabular}\;\times
\begin{tabular}{|c|}\hline
  \\ \hline
\end{tabular}\;=\;
\begin{tabular}{|c|}\hline
  \\ \hline
 \\ \hline
\end{tabular}\;+\;
\begin{tabular}{|c|c|}\hline
  & \\ \hline
\end{tabular}\;.
\]
Evidently, these are nothing but the well known Littlewood-Richardson relations
on the Schur functions.
\medskip

\nin{\bf Example \teknum} Take $\lambda=(4,2,1)$, $k=1$, $r_1=2$, $m_1=2$,
that is we add a border strip in the second and the third rows of $\lambda$.
The main identity takes the form
\[
s_{(4,2,1)}s_{(2,2)} = s_{(4,3,3)}s_{(1)}+ s_{(4,2)}s_{(2,2,1)}.
\]
In the graphic form it reads
\[
\begin{tabular}{|c|c|c|c|}\hline
  & & & \\ \hline
  & &\multicolumn{2}{c}{}\\ \cline{1-2}
  &\multicolumn{3}{c}{}\\ \cline{1-1}
\end{tabular}\;\times\;
\begin{tabular}{|c|c|}\hline
  & \\ \hline
  & \\ \hline
\end{tabular}\;=\;
\begin{tabular}{|c|c|c|c|}\hline
  & & & \\ \hline
  & & &\multicolumn{1}{c}{}\\ \cline{1-3}
 & & &\multicolumn{1}{c}{}\\ \cline{1-3}
\end{tabular}\;\times\;
\begin{tabular}{|c|}\hline
  \\ \hline
\end{tabular}\;+\;
\begin{tabular}{|c|c|c|c|}\hline
  & & & \\ \hline
  & &\multicolumn{2}{c}{}\\ \cline{1-2}
\end{tabular}\;\times\;
\begin{tabular}{|c|c|}\hline
  & \\ \hline
  & \\ \hline
  &\multicolumn{1}{c}{}\\ \cline{1-1}
\end{tabular}\;.
\]
Removing the first row or the first column gives rise to a pair of new identities
\[
s_{(2,1)}s_{(2,2)} = s_{(3,3)}s_{(1)} +s_{(2)}s_{(2,2,1)}
\]
and
\[
s_{(3,1)}s_{(2,2)} = s_{(3,2,2)}s_{(1)} +s_{(4,2)}s_{(1,1)}.
\]

Before proving the next corollary, we should introduce new notation.
With a Young diagram $\lambda=(\lambda_1,\dots,\lambda_n)$ we
associate a coordinate system with $x$ and $y$ axes directed as
shown in the picture below
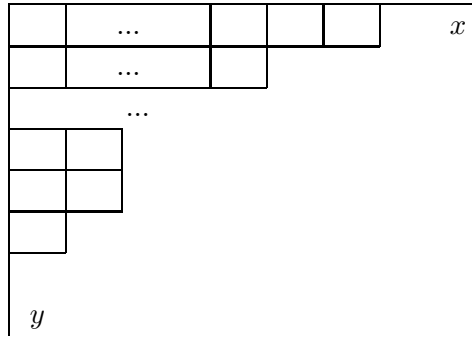
\[
\begin{array}{|c|c|c|c|c|c|c|cc}\hline
\hspace*{4mm} &\multicolumn{2}{c}{\hspace{5mm} . . .}&\hspace*{4mm} &\hspace*{4mm} & \hspace*{4mm}&\hspace*{4mm}&\hspace*{4mm}&x\\ \cline{1-7}
  & \multicolumn{2}{c}{\hspace{5mm} . . .}&\hspace*{4mm} & \multicolumn{1}{|c|}{} \\ \cline{1-5}
\multicolumn{5}{|c}{...}&\multicolumn{4}{c}{}\\ \cline{1-2}
&\hspace*{4mm} &\multicolumn{7}{c}{ }\\ \cline{1-2}
& &\multicolumn{7}{c}{}\\ \cline{1-2}
& \multicolumn{8}{c}{} \\ \cline{1-1}
\multicolumn{9}{|c}{}\\
\multicolumn{1}{|c}{y}&\multicolumn{8}{c}{}
\end{array}.
\]
The size of each box is accepted to be $1\times 1$, being measured
in the units of the $x$ and $y$ axes.

It is convenient to accept a different notation for components of a
given partition $\lambda$. Namely, we denote by $\xi_i$, where $1\le
i\le k\le n$, all {\it distinct} components of the partition
$\lambda$. That is $\lambda =
(\xi_1^{m_1},\xi_2^{m_2},\dots,\xi_k^{m_k})$ with some integers
$m_i\ge 1$, $m_1+\dots +m_k =n$. Note, that by definition
$\xi_i>\xi_j$ if $i<j$. Besides, it is convenient to set $\xi_{k+1}
= 0$. We also introduce a set of integers $y_i$, where $0\le i\le
k$, by the rule
\[
y_0=0,\quad y_i=m_i+y_{i-1},\quad 1\le i \le k.
\]
An {\it inner corner} of a diagram $\lambda$ is a point with
coordinates $(\xi_i,y_{i-1})$ with respect to the above coordinate
system. The collection of all inner corners will be called the {\it
inner corner set} of the diagram $\lambda$. So, the inner corner set
$\frak C_\lambda$ of the Young diagram
$\lambda=(\xi_1^{m_1},\dots,\xi_k^{m_k})$ consists of the following
$k+1$ points $\alpha_i$ \be {\frak C_\lambda}=\{\alpha_i =
(\xi_i,y_{i-1})\mid 1\le i\le k+1\}. \label{la-corn} \ee For
example, the inner corner set of the Young diagram $(6,5,2^2,1)$
includes five elements: $(6,0)$, $(5,1)$, $(2,2)$, $(1,4)$ and
$(0,5)$.

By the above definition, the inner corner set of any non-empty Young
diagram $\lambda$ is a non-empty set, containing at least two
elements --- the points $(\xi_1,0)$ and $(0,\ell(\lambda))$. Note
that knowing the inner corner set of a diagram allows one to restore
the diagram itself.

Introduce now the vertical and horizontal shifts of inner corners.
Let $\alpha_i=(\xi_i, y_{i-1})$, where $\xi_i\not=0$, be an inner
corner of a partition $\lambda=(\xi_1^{m_1},\dots,\xi_k^{m_k})$. The
{\it horizontal shift} $h^\pm_i$ of the corner $\alpha_i$ by $\pm 1$
means increasing or decreasing the component $\xi_i$ by $1$. If
$\xi_i+1=\xi_{i-1}$ or $\xi_i-1=\xi_{i+1}$, then the corresponding
rows of the diagram are united:
\[
\lambda=(\dots,\,\xi_{i-1}^{m_{i-1}},\xi_i^{m_i},\,\dots)
\stackrel{h_i^+}{\longrightarrow}
\left\{
\begin{array}{l@{\quad}c@{\quad}l}
(\dots,\,\xi_{i-1}^{m_{i-1}},(\xi_i+1)^{m_i},\,\dots)&{\rm if}&
\xi_{i-1}-\xi_i\ge 2\\
\rule{0pt}{7mm}
(\dots,\,\xi_{i-1}^{m_{i-1}+m_i},\,\dots)&{\rm if}&
\xi_{i-1}-\xi_i = 1,
\end{array}\right.
\]
\[
\lambda=(\dots,\,\xi_i^{m_i},\xi_{i+1}^{m_{i+1}},\,\dots)
\stackrel{h_i^-}{\longrightarrow}
\left\{
\begin{array}{l@{\quad}c@{\quad}l}
(\dots,\,(\xi_i-1)^{m_i},\xi_{i+1}^{m_{i+1}},\,\dots)&{\rm if}&
\xi_i-\xi_{i+1}\ge 2\\
\rule{0pt}{7mm}
(\dots,\,\xi_{i+1}^{m_i+m_{i+1}},\,\dots)&{\rm if}&
\xi_i-\xi_{i+1} = 1.
\end{array}\right.
\]
The other components of $\lambda$ preserve their values.

Similarly, the {\it vertical shift} $v_i^\pm$ of the corner
$\alpha_i=(\xi_i,y_{i-1})$, where $y_{i-1}\not=0$, by $\pm 1$
affects the exponents $m_i$ and $m_{i-1}$ in the following way
\[
(\dots,\,\xi_{i-1}^{m_{i-1}},\xi_i^{m_i},\,\dots)
\stackrel{v_i^-}{\longrightarrow}
\left\{
\begin{array}{l@{\quad}c@{\quad}l}
(\dots,\,\xi_{i-1}^{m_{i-1}-1},\xi_i^{m_i+1},\,\dots)&{\rm if }& m_{i-1}\ge 2\\
\rule{0pt}{7mm}
(\dots,\,\xi_{i-2}^{m_{i-2}},\xi_i^{m_i+1},\,\dots)&{\rm if }& m_{i-1} =1,
\end{array}\right.
\]
\[
(\dots,\,\xi_{i-1}^{m_{i-1}},\xi_i^{m_i},\,\dots)
\stackrel{v_i^+}{\longrightarrow}
\left\{
\begin{array}{l@{\quad}c@{\quad}l}
(\dots,\,\xi_{i-1}^{m_{i-1}+1},\xi_i^{m_i-1},\,\dots)&{\rm if }& m_i\ge 2\\
\rule{0pt}{7mm}
(\dots,\,\xi_{i-1}^{m_{i-1}+1},\xi_{i+1}^{m_{i+1}},\,\dots)&{\rm if }& m_i =1.
\end{array}\right.
\]
The other components of $\lambda$ remain unchanged.

Note that we do not define the horizontal shifts for the corner
$(0,\ell(\lambda))$ and vertical shifts for the corner $(\xi_1,0)$.

For example, for partition $\lambda=(6,5,2^2,1)$, the horizontal
shift of the corner $\alpha_3 =(2,2)$ by $+1$ and the vertical shift
of the corner $\alpha_2=(5,1)$ by $-1$ lead to the following
transformations:
\[
\lambda\;\stackrel{h_3^+}{\longrightarrow}(6,5,3^2,1),\qquad
\lambda\;\stackrel{v_2^{-}}{\longrightarrow} (5^2,2^2,1).
\]

Define now two transformations of any partition $\lambda$ generated
by shifts of the inner corners of the corresponding Young diagram.
\begin{definition}
\label{def:1} {\rm Let $\lambda$ be a partition and $\alpha_i =
(\xi_i,y_{i-1})$ an inner corner of the Young diagram $\lambda$.
Make the {\it horizontal} shift by $+1$ of all the inner corners
situated {\it above} $\alpha_i$ in the diagram $\lambda$ (that is
the corners $(\xi_{j},y_{j-1})$ with $j<i$). Besides, make the {\it
vertical} shift by $-1$ of all the inner corners situated {\it
below} $\alpha_i$ (that is, the corners $(\xi_j,y_{j-1})$ with
$j>i$). The corner $\alpha_i$ keeps its position unchanged. The
Young diagram thus obtained will be denoted $\lambda^+_-(\alpha_i)$.
In a similar way, shifting the corners above $\alpha_i$ by $-1$ in
the horizontal direction and those below $\alpha_i$ by $+1$ in the
vertical direction, we get the diagram $\lambda^-_+(\alpha_i)$. }
\end{definition}
Here is an example of the above procedures for the partition
$\lambda=(6,5,2^2,1)$ and the inner corner $\alpha_3=(2,2)$:
\[
\lambda=(6,5,2^2,1)\quad\Rightarrow\quad
\lambda^+_-(\alpha_3) = (7,6,2,1),\quad \lambda^-_+(\alpha_3) = (5,4,2^3,1).
\]

\begin{corollary}
\label{cor:f} {\rm Let $\lambda=(\xi_1^{m_1},\dots,\xi_k^{m_k})$ be
an arbitrary partition and let $\frak C_\lambda$ be the inner corner
set of the Young diagram $\lambda$. Then the following identity
holds true
\begin{equation}
s_\lambda s_\lambda =\sum_{\alpha\in {\frak C_\lambda}}s_{\lambda^+_-(\alpha)}
s_{\lambda^-_+(\alpha)}.
\label{sl-sl}
\end{equation}
This identity generalizes (\ref{br-kir}) to the case of an arbitrary partition.}
\end{corollary}
\noindent{\bf Proof.} Let $\lambda=(\xi_1^{m_1},\dots,\xi_k^{m_k})$
be an arbitrary partition of height $\ell(\lambda) = n$. We
introduce an auxiliary partition $\nu$ with $n+1$ components
\[
\nu=(\xi_1+1,\xi_1^{m_1},\xi_2^{m_2},\dots,\xi_k^{m_k}).
\]
On adding to the diagram $\nu$ all possible {\it strictly vertical} border strips,
we get the partition
\[
\nu^+=((\xi_1+1)^{m_1+1},(\xi_2+1)^{m_2},\dots,(\xi_k+1)^{m_k}).
\]
The inner corner sets of the new partitions are
\begin{eqnarray*}
&&{ \frak C_{\nu}} = (\xi_1+1,0)\cup \{(\xi_i, y_{i-1}+1)\,,\,1\le i\le k+1\}\\
\rule{0pt}{5mm}
&&{ \frak C_{\nu^+}} = \{(\xi_1+1,0),(0,y_k+1)\}\cup \{(\xi_i+1, y_{i-1}+1)\,,\,2\le i\le k\}
\end{eqnarray*}

Now we apply identity (\ref{main-id}) of Proposition \ref{prop:1} to
the product of the Schur functions $s_\nu s_{\nu^+\!\downarrow}$ and
then we use Corollary \ref{cor:2} in order to remove the first line
of length $\xi_1+1$ from the diagram $\nu$:
\[
\nu\mapsto \bar\nu = (\xi_1^{m_1},\dots,\xi_k^{m_k}) = \lambda.
\]
Besides, as follows from (\ref{srez}), $\nu^+\!\!\downarrow =
\lambda$. So, in our case, the left hand side of identity
(\ref{cor-id}) in Corollary \ref{cor:2} reads $s_{\bar\nu}
s_{\nu^+\!\downarrow} = s_\lambda s_\lambda$. We consider the right
hand side of (\ref{cor-id}) and verify that it coincides with that
of (\ref{sl-sl}).

The first term in the right hand side of (\ref{cor-id}) in our case
has the form $s_{\overline{\nu^+}}
s_{\raisebox{-1pt}{$_{\nu\downarrow}$}}$. Recall that the bar over
the symbol of partition means removing the first row of the
corresponding Young diagram. The inner corner sets of the diagrams
$\overline{\nu^+}$ and $\nu\!\!\downarrow$ are as follows
\begin{eqnarray*}
&&{ \frak C_{\overline{\nu^+}}} = \{(\xi_i+1, y_{i-1})\,,\,1\le i\le k\}\cup (0,y_k)\\
\rule{0pt}{5mm}
&&{ \frak C_{\nu\downarrow}} = \{(\xi_i-1, y_{i-1})\,,\,1\le i\le k\}\cup (0,y_k),
\end{eqnarray*}
and therefore, as follows from the structure of the inner corner set
$\frak C_\lambda$ (\ref{la-corn}) and Definition \ref{def:1},
\[
\overline{\nu^+} = \lambda^+_-(\alpha_{k+1}),\quad \nu\!\!\downarrow =
\lambda^-_+(\alpha_{k+1}),\qquad \alpha_{k+1} = (0,y_k).
\]

Consider now the sum over the partial peelings in (\ref{cor-id}). In
our case, this sum takes the form
\[
\sum_{p=1}^ks_{\overline{\nu^+}\downarrow^{(r_p)}}s_{\raisebox{-1pt}{$_{\nu\uparrow_{(r_p-1,1)}}$}}.
\]
The starting points $r_p$ of partial peelings in the diagram $\nu^+$
are the end points of the vertical border strips added to the
diagram $\nu$. The numbers $\{r_p\}$ are expressed in terms of
$\{y_p\}$ by the relation $r_p=y_{p-1}+2$ as illustrated in the
diagram below
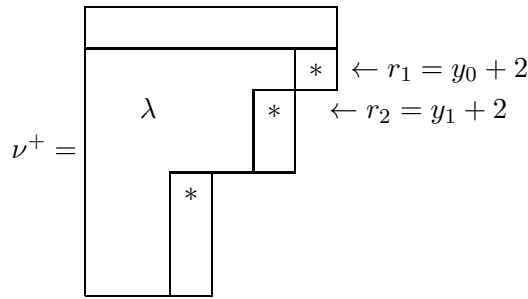
\[
\nu^+=
\begin{array}{|cccccc|cc}\cline{1-6}
\pha & \pha & \pha & \pha & \pha & \pha &\multicolumn{2}{c}{}\\ \cline{1-6}
\pha & \pha & \pha & \pha & \pha & \multicolumn{1}{|c|}{*}&\multicolumn{2}{c}{\leftarrow r_1=y_0+2}\\ \cline{5-6}
\pha & \lambda & \pha & \pha &\multicolumn{1}{|c|}{*}&\multicolumn{3}{c}{\leftarrow r_2=y_1+2}\\
\pha & \pha & \pha & \pha &\multicolumn{1}{|c|}{}&\multicolumn{3}{c}{}\\ \cline{3-5}
\pha & \pha &\multicolumn{1}{|c|}{*}&\multicolumn{5}{c}{}\\
\pha & \pha &\multicolumn{1}{|c|}{}&\multicolumn{5}{c}{}\\
\pha & \pha &\multicolumn{1}{|c|}{}&\multicolumn{5}{c}{}\\ \cline{1-3}
\end{array}
\]
Here the star signs mark the end points of the added border strips
--- the starting points $r_p$ of the partial down-peelings. As is
not difficult to see, the inner corner set of the diagram
$\overline{\nu^+}\!\!\downarrow^{(y_{p-1}+2)}$ has the following
structure
\[
{\frak C_{\overline{\nu^+}\!\downarrow^{(y_{p-1}+2)}}} =
\{(\xi_i+1,y_{i-1})\mid 1\le i\le p-1\}\cup (\xi_p,y_{p-1})\cup
\{(\xi_j,y_{j-1}-1)\mid p+1\le j \le k+1\}.
\]
By Definition \ref{def:1} this means that
\[
\overline{\nu^+}\!\!\downarrow^{(y_{p-1}+2)}
= \lambda^+_-(\alpha_p),\quad \alpha_p=(\xi_p,y_{p-1}).
\]
In analogous way we find that $\nu\!\!\uparrow_{(y_{p-1}+1,1)} =
\lambda^-_+(\alpha_p)$. Lastly, summation over $p$ gives the final
result (\ref{sl-sl}).\hfill\rule{6.5pt}{6.5pt}
\medskip

As an example we write down the bilinear relation for the square $s_{(3,2,1)}^2$:
\[
s_{(3,2,1)}s_{(3,2,1)} = s_{(4,3,2)}s_{(2,1)} + s_{(4,3)}s_{(2,1^3)} +s_{(4,1)}s_{(2^3,1)}
+s_{(3^2,2,1)}s_{(2,1)}.
\]

In what follows, we give a simple proof of the result (\ref{Dod})
\cite{FK}.

\begin{corollary}
\label{cor:3} {\rm {\bf $\cite{FK}$ }Let
$(\lambda_1,\lambda_2,\dots,\lambda_{n+1})$ be a partition with an
integer $n>0$. Then the following identity holds true \be
s_{(\lambda_2,\dots,\lambda_{n+1})}s_{(\lambda_1,\dots,\lambda_n)} =
s_{(\lambda_1+1,\dots,\lambda_n+1)}s_{(\lambda_2-1,\dots,\lambda_{n+1}-1)}
+s_{(\lambda_2,\dots,\lambda_n)}s_{(\lambda_1,\dots,\lambda_{n+1})}.
\label{ful-kl} \ee }
\end{corollary}

\nin{\bf Proof.} The result is based on identity (\ref{main-id}) and
the following steps.
\begin{enumerate}
\item Given a partition $\lambda=(\lambda_1,\lambda_2,\dots,\lambda_{n+1})$, we construct
an auxiliary partition
\[
\hat\lambda=(\lambda_1+1,\lambda_2,\dots,\lambda_{n+1})
\]
and take it as the initial partition for Proposition~\ref{prop:1}.
\item
Then we add to $\hat\lambda$ the connected
border strip from the second row till the last one ($k=1$, $r_1=2$, $m=n$) and get the
partition (see (\ref{la-pol}))
\[
\hat\lambda^+=(\lambda_1+1,\lambda_1+1,\lambda_2+1,\dots,\lambda_n+1).
\]
\item
Peeling the complete border strip off and partial peelings from the
end point of the added strip result in the following partitions (see
(\ref{srez}), (\ref{l-p-sr}) and (\ref{l-up-sr})):
\begin{eqnarray*}
&&\hat\lambda^+\!\!\downarrow = (\lambda_1,\lambda_2,\dots,\lambda_n)\\
&&\hat\lambda\!\!\downarrow = (\lambda_2-1,\dots,\lambda_{n+1}-1)\\
&&\hat\lambda^+\!\!\downarrow^{(2)} = (\lambda_1+1,\lambda_2,\dots, \lambda_n)\\
&&\hat\lambda\!\uparrow_{(1,\lambda_1-\lambda_2+1)} = (\lambda_1,\lambda_2,\dots,\lambda_{n+1}).
\end{eqnarray*}
\item
Lastly, the identity (\ref{main-id}) for the above Schur functions
gives
\begin{eqnarray*}
s_{(\lambda_1+1,\lambda_2,\dots,\lambda_{n+1})}
s_{(\lambda_1,\dots,\lambda_n)}&=&
s_{(\lambda_1+1,\lambda_1+1,\lambda_2+1,\dots,\lambda_n+1)}
s_{(\lambda_2-1,\dots,\lambda_{n+1}-1)}\\
& +&
s_{(\lambda_1+1,\lambda_2,\dots,\lambda_n)}
s_{(\lambda_1,\lambda_2,\dots,\lambda_{n+1})}.
\end{eqnarray*}
Removing from the above identity the first row $(\lambda_1+1)$ in
accordance with Corollary~\ref{cor:2}, we come to the result desired
(\ref{ful-kl}).\hfill\rule{6.5pt}{6.5pt}
\end{enumerate}

\noindent {\bf Note added in proof.} After this paper had been
accepted for publication, M.~Fulmek communicated to us that identity
(\ref{main-id}) can be proved in another way, as a corollary of
Lemma 16 in \cite{FK} (for details, see \cite{F}).

\section*{Acknowledgement}
This work was jointly supported by RFBR and
CNRS, grants 09-01-93107-NCNIL-a and GDRI-471. The work of P.P. and
P.S. was also supported by the RFBR grant 08-01-00392-a. 
The work of P.P. was supported by the grants of RFBR and SU-HSE
09-01-12185-ofi-m (09-09-0010) and by the SU-HSE grant 09-01-0026.

\end{document}